\documentclass[twoside,a4paper,12pt]{article}
\usepackage{amsmath,amsthm,amsfonts,amssymb}
\usepackage{mathrsfs}
\makeatletter
\def\@seccntformat#1{\csname the#1\endcsname.\quad}% To put a period
\renewcommand\section{\@startsection {section}{1}{\z@}%
                                   {-3.5ex \@plus -1ex \@minus -.2ex}% 
                                   {2.3ex \@plus.2ex}%
                                   {\normalfont\bf\center }}
\renewcommand\subsection{\@startsection {subsection}{1}{\z@}%
                                   {-3.5ex \@plus -1ex \@minus -.2ex}% 
                                   {2.3ex \@plus.2ex}%
                                   {\normalfont\bf}}
\makeatother

\date{\today}

\newcommand{\bR}{\mathbf{R}}

\newcommand{\bC}{\mathbf{C}}

\newcommand{\pa}{\partial}

\newcommand{\la}{\langle}
\newcommand{\ra}{\rangle}
\newcommand{\iti}{\boldsymbol{1}}
\newcommand{\sF}{\mathscr{F}}
\newtheoremstyle{new-thm}
 {3pt}
 {3pt}
 {\it}
 {0pt} %{\parindent}
 {\bf}
 {.}
 {.5em}
 {}
\newtheoremstyle{new-def}
 {3pt}
 {3pt}
 {\rm}
 {0pt} %{\parindent}
 {\bf}
 {.}
 {.5em}
 {}

\theoremstyle{new-thm}
  \newtheorem{thm}{Theorem}
  
  \newtheorem{lemma}[thm]{Lemma}
  \newtheorem{prop}[thm]{Proposition}
  
\theoremstyle{new-def}

  \newtheorem{rem}[thm]{Remark}

\numberwithin{equation}{section}
\numberwithin{thm}{section}

\begin{document}

\vspace*{1cm}
\begin{center} {\Large\bf 
Hitting times to spheres of Brownian motions\\
with and without drifts}
\end{center} 

\bigskip

\begin{center} Yuji Hamana and Hiroyuki Matsumoto \end{center} 

\bigskip

\begin{quote} {\bf Abstract.} 
Explicit formulae for the densities of the first hitting times 
to the sphere of Brownian motions with drifts are given.  
We need to consider the joint distributions of the first hitting 
times to the sphere and the hitting positions of 
the standard Brownian motion and explicit expression for 
their Laplace transforms are given, which are different from 
the known formulae in the literature and are of independnt interest. 

%%; Secondary 33C10(??), 44A10(??), \\
{\it keywords}\ : \ Brownian motion, first hitting time, 
Bessel process
\end{quote}

\section{Introduction}

For $d\geqq2$, we denote by $B=\{B_t\}_{t\geqq0}$ 
a standard $d$-dimensional Brownian motion 
starting from a fixed point $x\in\bR^d$ 
which is defined on a probability space 
$(\Omega,\sF,P)$.  
We throughout assume $x\ne0$.   
Letting $v\in\bR^d$ be a non-zero constant vector, 
we consider a Brownian motion $B^{(v)}=\{B^{(v)}_t\}_{t\geqq0}$ 
with drift $v$ given by $B_t^{(v)}=B_t+v t$.  
It is a very simple fundamental diffusion process, but 
we sometimes encounter difficulty to obtain explicit formulae on it. 

In this paper we consider the first hitting times $\sigma$ and 
$\sigma^{(v)}$ of $B$ and $B^{(v)}$, respectively, 
to the sphere $S^{d-1}_r$ with radius $r>0$ 
and centered at the origin.   
The main purpose is to give an explict expression for the density 
of $\sigma^{(v)}$.   
For this we need to give an explicit expression for the joint 
Laplace transform of the density of 
$(\sigma,B_\sigma)\in(0,\infty)\times S^{d-1}$.   
The formula for $(\sigma,B_\sigma)$ obtained in this article 
is of quite different form 
from the formulae obtained by Aizenman-Simon \cite{AS} and 
Wendel \cite{wendel}.  

The density $p_\nu(t;x)$, $\nu=\frac{d-2}{2}$ being the index, 
of $\sigma$ has been studied from old times.  
See \cite{K, PY} and the references therein 
for the Laplace forms and related topics.   
Recently, Byczkowski and Ryznar \cite{BR}, Uchiyama \cite{U} and 
the authors of the present paper \cite{HM-I, HM-T, HM-E} 
have studied the explicit expressions 
and the asymptotics of the densities themselves 
and the tail probabilities.   

The density for $\sigma^{(v)}$ is expressed in terms of 
the densities $p_\mu(t;x)$'s (of different dimensions).   
Moreover, using the previous results for $\sigma$, 
we show the asymptotics 
of the tail probabilities for $\sigma^{(v)}$.  

Our main results are the following.  

\begin{thm} \label{1t:main}
When $d=2,$ the density $p_0^{(v)}(t;x)$ for $\sigma^{(v)}$ 
is given by 
\begin{equation} \label{1e:main-1} \begin{split}
p_0^{(v)}(t;x) =  e^{-\la v,x \ra-\frac{1}{2}|v|^2t} 
 & \Bigl\{ I_0(|v|r) p_0(t;x) \\
 & + \sum_{n=1}^\infty n\; 
C_n^0\Bigl(\frac{\la v,x \ra}{|v|\cdot |x|}\Bigr) 
I_n(|v|r) \frac{|x|^n}{r^n}\; p_n(t;x) \Bigr\}.
\end{split} \end{equation}
When $d\geqq3,$ it is given by 
\begin{equation} \label{1e:main-2} \begin{split} 
p^{(v)}_\nu(t;x) = & 2^\nu \Gamma(\nu) 
 e^{-\la v,x \ra-\frac{1}{2}|v|^2t} \\
 & \times \sum_{n=0}^\infty (\nu+n) 
C_n^\nu\Bigl(\frac{\la v,x \ra}{|v|\;|x|}\Bigr) 
I_{\nu+n}(|v|r) \frac{|x|^n}{|v|^\nu r^{\nu+n}} 
\; p_{\nu+n}(t;x).  
\end{split} \end{equation}
Here $I_\mu\; (\mu\geqq0)$ is 
the modified Bessel function of the first kind 
and $C_n^\nu$ is the Gegenbauer polynomial{\rm .}  
\end{thm}

We refer to Magnus-Oberhettinger-Soni \cite{MOS} and 
Watson \cite{W} about the special functions.   
We note again that explicit expressions and the asymptotic behavior 
of $p_\nu(t;x)$ are known. 

For the asymptotic behavior of the tail probabilities, 
we show the following.  
To mention the result, we recall (\cite{BR, HM-I, U}) that 
\begin{equation*}
p_0(t;x)=\frac{L(0)}{t(\log t)^2} (1+o(1)) 
\quad \text{and} \quad 
p_\nu(t;x)=\frac{L(\nu)}{t^{\nu+1}}(1+o(1))\quad (\nu>0)
\end{equation*}
holds as $t\to\infty$, 
where $L(0)=2\log\frac{|x|}{r}$ and, for $\nu>0$, 
\begin{equation*}
L(\nu)= \frac{r^{2\nu}}{2^\nu \Gamma(\nu)} 
\Bigl( 1 - \Bigl(\frac{r}{|x|}\Bigr)^{2\nu}\Bigr).
\end{equation*}

\begin{thm} \label{1t:tail}
Assume $|x|>r.$  
Then{\rm ,} one has 
\begin{equation*}
P(t<\sigma^{(v)}<\infty) = \frac{2L(0)}{|v|^2} I_0(|v|r) 
e^{-\la v,x \ra} \frac{e^{-\frac{1}{2}|v|^2t}}{t(\log t)^2}
(1+o(1))
\end{equation*}
as $t\to\infty$ when $d=2,$ and 
\begin{equation*}
P(t<\sigma^{(v)}<\infty) = 
\frac{2^{\nu+1}L(\nu)\Gamma(\nu)}{|v|^2} 
\frac{I_\nu(|v|r)}{(|v|r)^{\nu}} e^{-\la v,x \ra} 
\frac{e^{-\frac{1}{2}|v|^2t}}{t^{\nu+1}}(1+o(1))
\end{equation*}
when $d\geqq3.$ 
\end{thm}

This paper is organized as follows.   
In the next section we give some estimates 
for the modified Bessel function and the Gegenbauer polynomial.  
In Section 3 we present an explicit form for the Laplace transform 
of $p^{(v)}_\nu(t;x)$ and, admitting it as proved, 
we give a proof of Theorem \ref{1t:main}.  
In Section 4 we show how to compute the Laplace form.   
In the last Section 5 we prove Theorem \ref{1t:tail}.  

We can apply the results in this paper 
to a study on the Wiener sausage of the Brownian motion 
with drift.  It will be discussed in a separate paper.  

\section{Preliminary estimates}

In this section we show some estimates 
for the modified Bessel function $I_\mu$ and 
for the Gegenbauer polynomial $C_n^\nu$.   

We firstly show an estimate for $I_\mu$: 
\begin{equation*}
I_\mu(\xi) = \sum_{m=0}^\infty 
\frac{(\xi/2)^{\mu+2m}}{\Gamma(m+1)\Gamma(m+\nu+1)}.
\end{equation*}

\begin{lemma} \label{2l:for-i}
For $\mu\geqq0$ and $n\geqq1,$ one has 
\begin{equation} \label{2e:for-i}
\xi^{-\mu}I_{\mu+n}(\xi) \leqq 
\frac{\xi^n}{2^{\mu+n}\Gamma(\mu+n+1)} e^\xi,\quad \xi>0.
\end{equation}
\end{lemma}

\noindent{\it Proof}.\quad
Note that $\Gamma(p+q)\geqq \Gamma(p+1) \Gamma(q)$ holds 
for $p\geqq0$ and $q\geqq1$, which can be seen from 
\begin{equation*}
\frac{\Gamma(p+1)\Gamma(q)}{\Gamma(p+q)} = 
p B(p,q) \leqq p \int_0^1 x^{p-1} dx = 1, \quad p>0.
\end{equation*}
Then we have 
\begin{align*}
\xi^{-\mu} I_{\mu+n}(\xi) & \leqq \frac{\xi^n}{2^{\mu+n}} 
\sum_{n=0}^\infty 
\frac{(\xi/2)^{2n}}{(\Gamma(m+1))^2\Gamma(\mu+n+1)} \\
 & \leqq \frac{\xi^n}{2^{\mu+n}\Gamma(\mu+n+1)} 
\Bigl( \sum_{m=0}^\infty \frac{(\xi/2)^m}{m!} \Bigr)^2,
\end{align*}
which shows \eqref{2e:for-i}. \qquad$\square$

\bigskip

Next we give an estimate for the Gegenbauer polynomial $C_n^\nu$.  
When $\nu>0$, it is given by 
\begin{equation*}
C_n^\nu(\xi) = \frac{1}{\Gamma(\nu)} \sum_{m=0}^{[n/2]} (-1)^m 
\frac{\Gamma(\nu+n-m)}{m! (n-2m)!} (2\xi)^{n-2m},
\end{equation*}
which is characterized by the relation 
\begin{equation*}
(1-2t\xi+t^2)^{-\nu} = \sum_{n=0}^\infty C_n^\nu(\xi) t^n.
\end{equation*}
When $\nu=0$, $C_0^0(\xi)=1$ and, when $n\geqq 1$, $C_n^0$ is given by 
\begin{equation*}
C_n^0(\xi)=\sum_{m=0}^{[n/2]} (-1)^m 
\frac{\Gamma(n-m)}{\Gamma(m+1)\Gamma(n-2m+1)} (2\xi)^{n-2m}.
\end{equation*}

\begin{lemma} \label{2l:for-g}
For $\alpha\in\bR$ with $|\alpha|\leqq 1,$ 
$\nu\geqq0$ and $n\geqq1,$ one has 
\begin{equation} \label{2e:for-g}
|C_n^\nu(\alpha)| \leqq \rho_\nu \frac{4^n \Gamma(\nu+n)}{n!},
\end{equation}
where $\rho_0=1$ and $\rho_\nu=(\Gamma(\nu))^{-1}$ for $\nu>0.$
\end{lemma}

\noindent{\it Proof}.\quad
When $\nu=0$, since $C_n^0(\cos\theta)=\frac{2}{n}\cos(n\theta)$, 
we have 
\begin{equation*}
|C_n^0(\alpha)|\leqq\frac{4^n}{n}=\rho_0\frac{4^n \Gamma(n)}{n!}.
\end{equation*}
When $\nu>0$, we have 
\begin{equation*}
\frac{|C_n^\nu(\alpha)|}{\Gamma(\nu+n)} \leqq 
\frac{1}{\Gamma(\nu)} \sum_{m=0}^{[n/2]} 
\frac{2^{n-2m} \Gamma(\nu+n-m)}{\Gamma(\nu+n)m! (n-2m)!}.
\end{equation*}
If $1\leqq m \leqq [n/2]$, it holds that 
\begin{equation*} \begin{split}
\frac{\Gamma(\nu+n-m)}{\Gamma(\nu+n) (n-2m)!} & \leqq 
\frac{1}{(n-1)(n-2)\cdots(n-m)\cdot(n-2m)!} \\
 & = \frac{(n-m-1)\cdots(n-(2m-1))}{(n-1)!} 
\leqq \frac{n^m}{n!}.
\end{split} \end{equation*}
Hence we get 
\begin{equation*} \begin{split}
\frac{|C_n^\nu(\alpha)|}{\Gamma(\nu+n)} & \leqq \frac{1}{\Gamma(\nu)}
\sum_{m=0}^{[n/2]} \frac{1}{m!} \frac{n^m}{n!} 2^{n-2m} \\ 
 & \leqq \frac{2^n}{\Gamma(\nu)n!} \sum_{m=0}^\infty 
\frac{1}{m!} \Bigl(\frac{n}{4}\Bigr)^m
\leqq \frac{2^n}{\Gamma(\nu)n!} e^{\frac{n}{2}} \leqq 
\frac{4^n}{\Gamma(\nu)n!}
\end{split} \end{equation*}
because $e^{\frac{1}{4}}\leqq2$.  \qquad$\square$

\section{Laplace transforms and proof of Theorem \ref{1t:main}}

We first reduce the computation for $\sigma^{(v)}$ 
to that for the joint distribution of $(\sigma, B_\sigma)$, 
the first hitting time to the sphere and the hitting position 
of the standard Brownian motion.   

By the Cameron-Martin theorem, we easily see 
\begin{equation*}
P(\sigma^{(v)} \leqq t) = e^{-\la v,x \ra-\frac12 |v|^2t} 
E[ e^{\la v, B_t \ra} \iti_{\{\sigma\leqq t\}}]
\end{equation*}
and
\begin{equation*}
E[e^{-\lambda \sigma^{(v)}}] = 
\lambda e^{-\la v,x \ra} 
\int_0^\infty e^{-(\lambda+\frac12 |v|^2)t} 
E[ e^{\la v, B_t \ra} \iti_{\{\sigma\leqq t\}}] dt, 
\end{equation*}
where $E$ denotes the expectation with respect to $P$.   
Moreover, letting $\sF_t=\sigma\{B_s, s\leqq t\}$, we see 
by the strong Markov property of Brownian motion
\begin{equation*}
E[ e^{\la v, B_t \ra} \iti_{\{\sigma\leqq t\}}] = 
E[ E[ e^{\la v, B_t \ra} | \sF_\sigma] \iti_{\{\sigma\leqq t\}}]
 = E[ e^{\la v, B_\sigma \ra + \frac12|v|^2(t-\sigma)} 
\iti_{\{\sigma\leqq t\}}]
\end{equation*}
and 
\begin{equation*}
E[e^{-\lambda\sigma^{(v)}}]=\lambda e^{-\la v,x \ra} 
\int_0^\infty e^{-\lambda t} E[ 
e^{\la v,B_\sigma \ra-\frac12 |v|^2\sigma} 
\iti_{\{\sigma\leqq t\}}] dt.
\end{equation*}
\indent
From this identity we obtain the following

\begin{prop} \label{2p:lap}
For any $\lambda>0$, one has 
\begin{equation} \label{3e:lap}
E[e^{-\lambda \sigma^{(v)}}] = 
e^{-\la v,x \ra} E[e^{\la v, B_\sigma \ra-
(\lambda+\frac12 |v|^2)\sigma}].
\end{equation}
\end{prop}

\noindent{\it Proof}.\quad
We have shown 
\begin{equation*}
E[e^{-\lambda \sigma^{(v)}}] = \lambda e^{-\la v,x \ra}
\int_0^\infty e^{-\lambda t}dt \int_0^t 
E[e^{\la v,B_\sigma \ra-\frac12 |v|^2\sigma} |\sigma=s] 
p_\nu(s;x) ds, 
\end{equation*}
where $p_\nu(s;x)$ is the density of $\sigma$ (see Sect. 1).   
Changing the order of integrations, we obtain 
\begin{equation*} \begin{split}
E[e^{-\lambda \sigma^{(v)}}] & = \lambda e^{-\la v,x \ra}
\int_0^\infty 
E[e^{\la v,B_\sigma \ra-\frac12 |v|^2\sigma} |\sigma=s] 
p_\nu(s;x) ds \int_s^\infty e^{-\lambda t}dt \\
 & = e^{-\la v,x \ra} E[ e^{\la v,B_\sigma\ra-
(\lambda+\frac12 |v|^2)\sigma}]. \qquad \qquad \square
\end{split} \end{equation*}

\medskip

For the right hand side of \eqref{3e:lap}, 
we show the following explicit expression.   
Denoting by $K_\mu$ the modified Bessel function of the second kind 
(the Macdonald function), 
we define the function $Z_\mu^{(v),\lambda}\ (\mu\geqq 0)$ by 
\begin{align*}
 & Z_\mu^{(v),\lambda}(\xi,\eta) = 
\frac{K_\mu(\xi\sqrt{2\lambda+|v|^2})}
{K_\mu(\eta\sqrt{2\lambda+|v|^2})}\qquad \text{if\ $\xi>\eta>0$}\\
\intertext{and}
 & Z_\mu^{(v),\lambda}(\xi,\eta) = 
\frac{I_\mu(\xi\sqrt{2\lambda+|v|^2})}
{I_\mu(\eta\sqrt{2\lambda+|v|^2})}\qquad \text{if\ $\eta>\xi>0$}.
\end{align*}
Since $K_\mu$ is decreasing and $I_\mu$ is increasing 
on $(0,\infty)$, 
$Z_\mu^{(v),\lambda}\leqq1$.  

\begin{prop} \label{3p:joint}
Let $\lambda>0.$    When $d=2,$ one has 
\begin{equation*}
E[e^{\la v,B_\sigma \ra-\lambda\sigma}] = 
I_0(|v|r) Z_0^{(0),\lambda}(|x|,r) + 
\sum_{n=1}^\infty n\; C_n^0(\alpha) I_n(|v|r) 
Z_n^{(0),\lambda}(|x|,r),
\end{equation*}
where $\alpha=\frac{\la v,x \ra}{|v|\cdot|x|}.$   
When $d\geqq3,$ one has 
\begin{equation*}
E[e^{\la v,B_\sigma \ra-\lambda\sigma}] = 
2^\nu \Gamma(\nu) \sum_{n=0}^\infty (\nu+n) C_n^\nu(\alpha) 
\frac{I_{\nu+n}(|v|r)}{(|v|\cdot|x|)^\nu} 
Z_{\nu+n}^{(0),\lambda}(|x|,r).
\end{equation*}
\end{prop}

Combining this proposition with \eqref{3e:lap}, we obtain 
\begin{align*} 
 & E[e^{-\lambda\sigma^{(v)}}] = \int_0^\infty e^{-\lambda t}
p^{(v)}(t;x) dt \\ 
=\ & e^{-\la v,x \ra} \Bigl\{ I_0(|v|r) 
Z_0^{(v),\lambda}(|x|,r) + \sum_{n=1}^\infty n\; C_n^0(\alpha)
I_n(|v|r) Z_n^{(v),\lambda}(|x|,r) \Bigr\} \\
\intertext{when $d=2,$ and} 
 & E[e^{-\lambda\sigma^{(v)}}] = 
e^{-\la v,x \ra} 2^\nu  \Gamma(\nu) \sum_{n=0}^\infty (\nu+n)
C_n^\nu(\alpha) \frac{I_{\nu+n}(|v|r)}{(|v|\cdot|x|)^\nu} 
Z_{\nu+n}^{(v),\lambda}(|x|,r).
\end{align*}
when $d\geqq3$.  

We postpone a proof of Propositon \ref{3p:joint} 
to the next section and give a proof of Theorem \ref{1t:main}.  
It is well known (cf. \cite{K}) that, 
the density $p_\mu(s;x)$ of $\sigma_\mu$, 
the first hitting time of a Bessel process with index $\mu$ 
starting from $|x|$, is characterized by 
\begin{equation} \label{3e:lap-known} 
E[e^{-\lambda \sigma_\mu}] = \int_0^\infty e^{-\lambda s} 
p_\mu(s;x) ds = 
\frac{r^{\mu}}{|x|^\mu} Z_\mu^{(0),\lambda}(|x|,r).
\end{equation}
We use the same notation for the density 
since our main concern is on the special case 
where the index $\mu$ is a half integer, 
and there is no fear of confusion.  

To prove Theorem \ref{1t:main}, 
we compute the Laplace transform of the right hand side 
of \eqref{1e:main-1}, \eqref{1e:main-2} 
by changing the order of the integrations and the infinite sums.  

Using the estimates given in the previous section, 
we have for $\nu\geqq0$ 
\begin{equation} \label{3e:joint} \begin{split}
\sum_{n=1}^\infty & (\nu+n) |C_n^\nu(\alpha)| 
 I_{\nu+n}(|v|r) \frac{|x|^n}{|v|^\nu r^{\nu+n}} 
\int_0^\infty e^{-(\lambda+\frac{1}{2}|v|^2)t} 
p_{\nu+n}(t;x) dt \\
 & \leqq \frac{\rho_\nu r^\nu}{2^\nu |x|^\nu} e^{|v|r} 
\sum_{n=1}^\infty \frac{(\nu+n) \Gamma(\nu+n) (2|v|r)^n}
{n! \Gamma(\nu+n+1)} Z_{\nu+n}^{(v),\lambda}(|x|,r).
\end{split} \end{equation}
Since $Z_\mu^{(v),\lambda}\leqq1$, the above is bounded by 
\begin{equation*}
\frac{\rho_\nu r^\nu}{2^\nu|x|^\nu} e^{|v|r} \sum_{n=0}^\infty 
\frac{(2|v|r)^n}{n!} = 
\frac{\rho_\nu r^\nu}{2^\nu|x|^\nu} e^{3|v|r}.
\end{equation*}
\indent
Hence, we may apply Fubini's theorem and 
see, from \eqref{3e:lap-known}, that 
the Laplace transforms of the right hand sides of 
\eqref{1e:main-1} and \eqref{1e:main-2} are equal to 
those of $p^{(v)}_\nu(t;x)$ in both cases.   

We have now shown Theorem \ref{1t:main}, 
admitting Proposition \ref{3p:joint} as proved.  

\section{Proof of Proposition \ref{3p:joint}}

In order to prove Proposition \ref{3p:joint}, 
we use the skew-product representation of Brownian motions.   
Let $R=\{R_t\}_{t\geqq0}$ be a $d$-dimensional Bessel 
process (with index $\nu=\frac{d-2}{2}$) and 
$\theta=\{\theta_t\}_{t\geqq0}$ be a Brownian motion 
on the unit sphere $S^{d-1}=S_1^{d-1}$ 
with $\theta_0=\frac{x}{|x|}$, 
and assume that $R$ and $\theta$ are independent.  
Recall that, embedding $S^{d-1}$ in $\bR^d$, 
we can realize $\theta$ as a solution 
of a stochastic differential equation, 
which is so-called Stroock's representation 
of a spherical Brownian motion.   

Set $S_t=\int_0^t(R_s)^{-2}ds$.   
Then, $\{R_t\theta_{S_t}\}_{t\geqq0}$ is a $d$-dimensional 
Brownian motion.   Hence, we have 
\begin{equation*}
E[e^{-\lambda\sigma+\la v,B_\sigma \ra}] = 
\int_0^\infty \int_0^\infty e^{-\lambda t} 
E_{\frac{x}{|x|}}^\theta[e^{r\la v,\theta_u \ra}] 
P_{\nu,|x|}(\tau\in dt, S_\tau\in du),
\end{equation*}
where $E_{\theta_0}^\theta$ denotes the expectation 
with respect to the probability law of $\theta$ 
starting from $\theta_0$, 
$P_{\nu,|x|}$ is the probability law of $\{R_t\}$ and 
$\tau$ is the first hitting time to $r$ of $\{R_t\}$.   

It is known (cf. \cite{BS} p.407) that
\begin{align*}
 & E_{\nu,|x|}[ e^{-\alpha\tau-\frac12 \beta^2S_\tau}] = 
\frac{|x|^{-\nu}K_{\sqrt{\nu^2+\beta^2}}(|x|\sqrt{2\alpha})}
{r^{-\nu}K_{\sqrt{\nu^2+\beta^2}}(r\sqrt{2\alpha})}
\qquad \text{\rm if}\quad |x|>r \\
\intertext{and}
 & E_{\nu,|x|}[ e^{-\alpha\tau-\frac12 \beta^2S_\tau}] = 
\frac{|x|^{-\nu}I_{\sqrt{\nu^2+\beta^2}}(|x|\sqrt{2\alpha})}
{r^{-\nu}I_{\sqrt{\nu^2+\beta^2}}(r\sqrt{2\alpha})}
\qquad \text{\rm if}\quad |x|<r , 
\end{align*}
where $E_{\nu,|x|}$ is the expectatation 
with respect to $P_{\nu,|x|}$.   

We obtain Proposition \ref{3p:joint} 
if we show the following.  
We can justify the change of order of the integration and 
the infinite sum by the same way as \eqref{3e:joint}. 

\begin{prop} \label{4p:sphere} 
Let $\xi>0$.   
Then{\rm ,} when $d=2,$ one has
\begin{equation} \label{4e:sphere-1}
E_{\theta_0}^\theta[ e^{\xi\la v,\theta_t\ra}] = 
I_0(|v|r) + \sum_{n=1}^\infty n\; C_n^0(\alpha)
e^{-\frac{1}{2}n^2t} I_n(|v|\xi)
\end{equation}
and{\rm ,} when $d\geqq3,$ 
\begin{equation} \label{4e:sphere-2}
E^\theta_{\theta_0}[e^{\xi\la v,\theta_t\ra}]=2^\nu \Gamma(\nu)
\sum_{n=0}^\infty (\nu+n)
C_n^\nu(\alpha) e^{-\frac12 n(n+2\nu)t} 
\frac{I_{\nu+n}(|v|\xi)}{(|v|\xi)^{\nu}}.
\end{equation}
\end{prop}

We see from this proposition that 
the Gegenbauer polynomial comes into our story 
through the following formula (cf. \cite[p.227]{MOS}): 
for $\alpha\in\bR,\xi>0,\mu>0$, 
\begin{equation} \label{4e:gegen}
e^{\alpha\xi}=2^\mu \Gamma(\mu) \sum_{n=0}^\infty 
(\mu+n) C_n^\mu(\alpha) \xi^{-\mu} I_{\mu+n}(\xi).
\end{equation}
\indent We first show that 
$f_\nu(t,\xi)=E^\theta_{\theta_0}[e^{\xi\la v,\theta_t\ra}]$
satisfies 
\begin{equation} \label{4e:heat}
\frac{\pa f_\nu}{\pa t}=-\frac12 \xi^2 \frac{\pa^2f_\nu}{\pa\xi^2}
-\frac{d-1}{2}\xi\frac{\pa f_\nu}{\pa\xi}+
\frac12 |v|^2\xi^2 f_\nu, \quad 
t>0,\ \xi>0,
\end{equation}
together with the boundary conditions
\begin{equation} \label{4e:bdry}
f_\nu(0,\xi)=e^{\xi\la v,\theta_0 \ra}, \quad
f_\nu(t,0)=1, \quad
\frac{\pa f_\nu}{\pa\xi}(t,0)=\la v,\theta_0 \ra 
e^{-\frac{d-1}{2}t}.
\end{equation}
\indent
For this purpose, we recall Stroock's representation 
of sperical Brownian motion (cf. \cite{S}).   
$\theta$ may be realized as a solution of 
the stochastic differential equation 
based on a $d$-dimensional Brownian motion 
$\{w_s=(w_s^1,w_s^2,...,w_s^d)\}_{s\geqq0}$ which is given by 
\begin{equation*}
d\theta_s^i = \sum_{j=1}^d 
(\delta_{ij}-\theta_s^i \theta_s^j) \circ dw_s^j, 
\qquad i=1,2,...,d.
\end{equation*}
Then, from a strihgtforward computation by using It\^o's formula, 
we can show \eqref{4e:heat}.  
It is easy to see \eqref{4e:bdry}.   

For simplicity we set $\beta=d-1$ 
and consider the function $g_\nu$ given by 
\begin{equation*}
g_\nu(t,\xi)=f_\nu\Bigl(2t,\frac{\xi}{|v|}\Bigr).
\end{equation*}
Then $g_\nu$ is a smooth function which satisfies 
\begin{align}
 & \frac{\pa g_\nu}{\pa t}=-\xi^2 \frac{\pa^2g_\nu}{\pa \xi^2}
-\beta \xi \frac{\pa g_\nu}{\pa\xi} + \xi^2g_\nu, 
\qquad t>0,\ \xi>0, \label{4e:diff-eq}
\intertext{and} 
 & g_\nu(0,\xi)=e^{\alpha\xi},\qquad g_\nu(t,0)=1, \qquad 
\frac{\pa g_\nu}{\pa\xi}(t,0)=
\alpha e^{-\beta t}.\label{4e:bdry-cond}
\end{align}
\indent 
If $u(t,\xi)=e^{-\lambda t}\phi(\xi)$ satisfies \eqref{4e:diff-eq}, 
we should have 
\begin{equation*}
\xi^2\phi''(\xi) + \beta\xi\phi'(\xi)-(\xi^2+\lambda)\phi(\xi)=0.
\end{equation*}
The system of the fundamental solutions 
of this second order differential equation is given by 
$\xi^{-\nu}I_{\sqrt{\lambda+\nu^2}}(\xi)$ and 
$\xi^{-\nu}K_{\sqrt{\lambda+\nu^2}}(\xi)$, 
where $\nu=\frac{\beta-1}{2}=\frac{d-2}{2}.$  
For the function $\phi$ to be smooth at $\xi=0$, 
we should choose $\xi^{-\nu}I_{\sqrt{\lambda+\nu^2}}(\xi)$.  
Moreover, $n=\sqrt{\lambda+\nu^2}-\nu$ should 
be a non-negative integer and $\lambda=n(n+2\nu)$.   

The following lemma is easily shown and we omit the proof.  

\begin{lemma} \label{4l:bessel}
{\rm (1)}\ 
The function 
$\varphi_{\nu,n}(\xi)=\xi^{-\nu}I_{\nu+n}(\xi)$ satisfies 
\begin{equation*}
\xi^2\varphi_{\nu,n}''(\xi)+\beta \xi \varphi_{\nu,n}'(\xi) - 
\xi^2 \varphi_{\nu,n}(\xi) = n(n+2\nu) 
\varphi_{\nu,n}(\xi), \qquad \xi>0.
\end{equation*}
{\rm (2)}\ One has 
\begin{equation*}
\varphi_{\nu,1}'(0)=\frac{1}{2^{\nu+1}\Gamma(\nu+2)} 
\qquad \text{and} \qquad
\varphi_{\nu,n}'(0)=0 \quad (n\ne1).
\end{equation*}
\end{lemma}

The following proposition 
immediately implies Proposition \ref{4p:sphere} 

\begin{prop} \label{4p:expl-rep}
When $d=2,$ one has 
\begin{equation} \label{4e:expl-rep1}
g_0(t,\xi)=I_0(\xi) + \sum_{n=1}^\infty n\; C_n^0(\alpha)
e^{-n^2t} I_n(\xi),\quad t\geqq0, \ \xi\geqq0
\end{equation}
and{\rm ,} when $d\geqq3,$ 
\begin{equation} \label{4e:expl-rep2} 
g_\nu(t,\xi)= 2^\nu \Gamma(\nu) \sum_{n=0}^\infty 
(\nu+n) C_n^\nu(\alpha) e^{-n(n+2\nu)t} \xi^{-\nu} I_{\nu+n}(\xi).
\end{equation}
\end{prop}

\noindent{\it Proof}.\quad
First of all we note that 
the sum on the right hand sides 
of \eqref{4e:expl-rep1} and \eqref{4e:expl-rep2} are 
absolutely convergent at each $(t,\xi)$, 
which is seen from \eqref{2e:for-i} and \eqref{2e:for-g} 
in a similar way to \eqref{3e:joint}.   

Letting $\varphi_{\nu,n}$ be the function 
defined in Lemma \ref{4l:bessel}, we set 
\begin{equation*}
h_{\nu}(t;\xi) = \sum_{n=1}^\infty (\nu+n) C_n^\nu(\alpha)
e^{-n(n+2\nu)t} \varphi_{\nu,n}(\xi).
\end{equation*}
We have shown that the sum on the right hand side is 
absolutely convergent.   
Moreover, noting 
\begin{equation*}
\varphi'_{\nu,n}(\xi)=\frac{n}{\xi}\varphi_{\nu,n}(\xi)
+\varphi_{\nu,n+1}(\xi),
\end{equation*}
we see, in a similar way to \eqref{3e:joint}, that
\begin{equation*}
\sum_{n=1}^\infty (\nu+n) C_n^\nu(\alpha)
e^{-n(n+2\nu)t} \varphi_{\nu,n}'(\xi)
\ \ \text{\rm and} \ \ 
\sum_{n=1}^\infty (\nu+n) C_n^\nu(\alpha)
e^{-n(n+2\nu)t} \varphi_{\nu,n}''(\xi)
\end{equation*}
converge uniformly on compact sets in $\xi\in(0,\infty)$ 
and are equal to $\frac{\pa}{\pa\xi}h_\nu(t,\xi)$ and 
$\frac{\pa^2}{\pa\xi^2}h_\nu(t,\xi)$, respectively.   

Next we look at $h_\nu(t,\xi)$ as a function in $t>0$.   
By \eqref{2e:for-i} and \eqref{2e:for-g} we have 
\begin{align*}
 & \sum_{n=1}^\infty (\nu+n) |C_n^\nu(\alpha)| n(n+2\nu)
e^{-n(n+2\nu)t} \varphi_{\nu,n}(\xi) \\
\leqq\ & \rho_\nu \sum_{n=1}^\infty (\nu+n) 
\frac{4^n\Gamma(\nu+n)}{n!} n(n+2\nu) 
\frac{\xi^n}{2^{\nu+n}\Gamma(\nu+n+1)} e^\xi \\
=\ & \frac{\rho_\nu}{2^\nu} e^\xi \Bigl\{ 
\sum_{n=1}^\infty \frac{(n-1)(2\xi)^n}{(n-1)!} + 
\sum_{n=1}^\infty \frac{(2\nu+1)(2\xi)^n}{(n-1)!} \Bigr\} \\
=\ & \frac{\rho_\nu \xi^2}{2^{\nu-2}} e^{3\xi} + 
\frac{\rho_\nu(2\nu+1)\xi}{2^{\nu-1}} e^{3\xi}
\end{align*}
and we may differetiate term by term to obtain 
\begin{equation*}
\frac{\pa}{\pa t}h_\nu(t,\xi) = 
-\sum_{n=1}^\infty (\nu+n) C_n^\nu(\alpha) e^{-n(n+2\nu)t} 
n(n+2\nu) \varphi_{\nu,n}(\xi).
\end{equation*}
\indent
Combining the identities above, 
we see that the function $g_\nu(t,\xi)$ given 
by \eqref{4e:expl-rep1} and \eqref{4e:expl-rep2} 
satisfies \eqref{4e:diff-eq}.   

The boundary condition \eqref{4e:bdry-cond} 
in the case of $d\geqq3$ may be checked by \eqref{4e:gegen} 
and the fact, $C_0^\nu(\alpha)=1$ and 
$C_1^\nu(\alpha)=2\nu\alpha$\ (cf. \cite[p.218]{MOS}).  

For a check when $d=2$, we rewrite \eqref{4e:gegen} as 
\begin{align*}
e^{\alpha\xi} & = 2^\mu \Gamma(\mu+1) \xi^{-\mu} I_\mu(\xi) + 
2^\mu \sum_{n=1}^\infty \Gamma(\mu) (\mu+n) C_n^\mu(\alpha) 
\xi^{-\mu} I_{\mu+n}(\xi) \\
 & = 2^\mu \Gamma(\mu+1) \Bigl\{ \xi^{-\mu} I_\mu(\xi) + 
\sum_{n=1}^\infty (\mu+n) \frac{C_n^\mu(\alpha)}{\mu} 
\xi^{-\mu} I_{\mu+n}(\xi) \Bigr\},
\end{align*}
which holds for any $\mu>0$.   
Note $\frac{C_n^\mu(\alpha)}{\mu}\to C_n^0(\alpha)$ 
as $\mu\downarrow0$.   
Then, by using \eqref{2e:for-i} and \eqref{2e:for-g}, 
we can show that we may apply the dominated convergence theorem 
and obtain 
\begin{equation*}
e^{\alpha\xi} = I_0(\xi) + \sum_{n=1}^\infty n\; C_n^0(\alpha) 
I_n(\xi),
\end{equation*}
which is exactly $g_0(0,\xi)=e^{\alpha\xi}.$   

Another boundary condition $g_0(t,0)=1$ follows from $I_0(0)=1$ 
and the other one $\frac{\pa g_\nu}{\pa\xi}(t,0)=\alpha e^{-t}$ 
does from the formula $C_1^0(\alpha)=\alpha$.  

We have now completed the proof of Proposition \ref{4p:expl-rep}.

\begin{rem}
The function $z^{-\nu}I_{\nu+n}(z)$ may be regarded as 
a holomorphic function on $\bC$.  
Hence, by using \eqref{2e:for-i} and \eqref{2e:for-g}, 
we can show that the functions 
\begin{equation*}
E[e^{z\la v,\theta_t \ra/|v|}] \quad \text{and} \quad 
\sum_{n=1}^\infty (\nu+n) C_n^\nu(z) e^{-\frac{1}{2}n(n+2\nu)t} 
z^{-\nu} I_{\nu+z}(z)
\end{equation*}
are holomorphic in $z\in\bC$.   
From this we obtain the Fourier-Laplace transform of 
the joint distribution of $(\sigma,B_\sigma)$.   
For example, when $d\geqq3$, we can show 
\begin{equation*}
E[e^{i\la v,B_\sigma\ra - \lambda\sigma}] = 
2^\nu \Gamma(\nu) \sum_{n=0}^\infty i^n(\nu+n) 
C_n^\nu(\alpha) \frac{J_{\nu+n}(|v|r)}{(|v|\cdot|x|)^\nu}
Z^{(0),\lambda}_{\nu+n}(|x|,r),
\end{equation*}
where $J_\mu$ is the usual Bessel function. 
\end{rem}

\section{Proof of Theorem \ref{1t:tail}}

We set 
\begin{equation*}
f^{(v)}(t;x)=P(t<\sigma^{(v)}<\infty)=
\int_t^\infty p^{(v)}(s;x)ds.
\end{equation*}
In order to apply Theorem \ref{1t:main}, 
we need to change the order of the integration 
and the summation.   

For this purpose, 
we note from \eqref{2e:for-i} and \eqref{2e:for-g} 
\begin{align*}
 & \sum_{n=1}^\infty (\nu+n) |C_n^\nu(\alpha)| 
\frac{I_{\nu+n}(|v|r)|x|^n}{|v|^\nu r^{\nu+n}}
\int_t^\infty e^{-\frac12 |v|^2s} p_{\nu+n}(s;x)ds \\
\leqq\ & \sum_{n=1}^\infty (\nu+n) \rho_\nu 
\frac{4^n \Gamma(\nu+n)}{n!} 
\frac{(|v|r)^ne^{|v|r}}{2^{\nu+n}\Gamma(\nu+n+1)}
\frac{|x|^n}{r^n} e^{-\frac{1}{2}|v|^2t} \\
\leqq\ & \frac{\rho_\nu}{2^\nu} e^{|v|r} \sum_{n=1}^\infty 
\frac{(2|v|\cdot|x|)^n}{n!} 
\leqq \frac{\rho_\nu}{2^\nu} e^{3|v|r}.
\end{align*}
Then we can apply Fubini's theorem and 
obtain from Theorem \ref{1t:main} 
\begin{equation*}
f^{(v)}(t;x)=e^{-\la v,x\ra}I_0(|v|r) \int_t^\infty 
e^{-\frac{1}{2}|v|^2s} p_0(s;x)ds + f_0(t)
\end{equation*}
when $d=2$, and when $d\geqq3$
\begin{equation*}
f^{(v)}(t,x)=e^{-\la v,x \ra} 2^\nu \Gamma(\nu+1) 
\frac{I_\nu(|v|r)}{|v|^\nu r^\nu} \int_t^\infty 
e^{-\frac{1}{2}|v|^2s} p_\nu(s;x)ds + f_\nu(t),
\end{equation*}
where the second terms on the right hand sides are given by 
\begin{equation*} \begin{split}
f_\nu(t) = & e^{-\la v,x \ra} 2^\nu \Gamma(\nu+1) \\
 & \times \sum_{n=1}^\infty (\nu+n) C_n^\nu(\alpha) 
\frac{|x|^n I_{\nu+n}(|v|r)}{|v|^\nu r^{\nu+n}} 
\int_t^\infty e^{-\frac{1}{2}|v|^2s}p_{\nu+n}(s;x)ds.
\end{split} \end{equation*}
\indent 
At first we prove the theorem when $\nu>0$.   
We have
\begin{equation*}
p_\nu(t;x) = \frac{L(\nu)}{t^{\nu+1}} (1+o(1))
\end{equation*}
and, by L'Hospital's rule, 
\begin{equation} \label{4e:tail} 
\int_t^\infty e^{-\frac{1}{2}|v|^2s} p_\nu(s;x)ds = 
\frac{2L(\nu)}{|v|^2} t^{-\nu-1} e^{-\frac{1}{2}|v|^2t}
(1+o(1)).
\end{equation}
Note that this identity holds for any $\nu>0$.   

In order to show that $f_\nu(t)$ is negligible when $\nu>0$, 
we use the argument in Sect.2 of \cite{HM-E}.
Then we obtain 
\begin{equation*} \begin{split}
\int_t^\infty e^{-\frac{1}{2}|v|^2s}p_{\nu+n}(s;x)ds & \leqq 
e^{-\frac{1}{2}|v|^2t} P(t<\sigma^{(\nu+n)}<\infty) \\
 & \leqq e^{-\frac{1}{2}|v|^2t} E_{\nu+n,|x|}[(R_t)^{-2(\nu+n)}],
\end{split} \end{equation*}
where $E_{\nu+n,|x|}$ denotes the expectation 
with respect to the probability law of the Bessel process 
$\{R_t\}_{t\geqq0}$ with index $\nu+n$ and starting from $|x|$.   

Using the explicit expression of the transition density 
of the Bessel process, we obtain 
\begin{equation*} \begin{split}
E_{\nu,|x|}[(R_t)^{-2(\nu+n)}] & = \frac{1}{(2t)^{\nu+n}} 
e^{-\frac{|x|^2}{2t}} \sum_{m=0}^\infty 
\frac{|x|^m}{\Gamma(\nu+n+m+1)(2t)^m} \\ 
 & \leqq \frac{1}{(2t)^{\nu+n}} e^{-\frac{|x|^2}{2t}} 
\sum_{m=0}^\infty 
\frac{|x|^m}{\Gamma(\nu+n+1)\Gamma(m+1)(2t)^m} \\
 & = \frac{1}{\Gamma(\nu+n+1)(2t)^{\nu+n}}.
\end{split} \end{equation*}
Hence, by using \eqref{2e:for-i} and \eqref{2e:for-g} again, 
we obtain for $n\geqq1$ and $t\geqq1$
\begin{equation*} \begin{split} 
 & t^{\nu+1} e^{\frac12 |v|^2t} (\nu+n) |C_n^\nu(\alpha)| 
\frac{|x|^n I_{\nu+n}(|v|r)}{|v|^\nu r^{\nu+n}} 
\int_t^\infty e^{-\frac{1}{2}|v|^2s} p_{\nu+n}(s;x)ds \\
 & \leqq \frac{\rho_\nu}{4^\nu} e^{|v|r} 
\frac{(|v|\cdot|x|)^n}{n! \Gamma(\nu+n+1)}.
\end{split} \end{equation*}
The quantity on the right hand side is independent of $t\geqq1$ 
and is summable in $n$.   
Therefore, since \eqref{4e:tail} holds when we replace 
$\nu$ by $\nu+n$, we can apply the dominated convergence theorem 
and obtain 
\begin{equation*}
\lim_{t\to\infty} t^{\nu+1}e^{\frac12 |v|^2t} f_\nu(t)=0.
\end{equation*}
\indent 
When $d=2$, we have 
\begin{equation*}
p_0(t;x)=\frac{L(0)}{t(\log t)^2}(1+o(1)), \qquad 
L(0)=2\log\frac{|x|}{r},
\end{equation*}
and 
\begin{equation*}
\int_t^\infty e^{-\frac{1}{2}|v|^2s} p_0(s;x)ds = 
\frac{L(0)}{t(\log t)^2} e^{-\frac{1}{2}|v|^2t} (1+o(1)).
\end{equation*}
By the same way as in the case of $\nu>0$, 
we can show by the dominated convergence theorem 
\begin{equation*}
\lim_{t\to\infty} t(\log t)^2 e^{\frac{1}{2}|v|^2t}f_0(t)=0
\end{equation*}
and we obtain the assertion of Theorem \ref{1t:tail} 
also in this case. 

\begin{rem}
The estimate for the tail probability $P(t<\sigma<\infty)$ 
has been firstly given in Byczkowski and Ryznar \cite{BR}.  
We need an explicit upper bound here.
\end{rem}

\section*{Acknowledgements}
The authors thank Professor Tatsuo Iguchi 
for his valuable suggestions.


\begin{thebibliography}{99}

\bibitem{AS} M.~Aizenman and B.~Simon, 
Brownian motion and Harnack inequality 
for Schr\"odinger operators, 
Commun. Pure Appl. Math., {\bf 35} (1982), 209--273.  

\bibitem{BS} A.~N.~Borodin and P.~Salminen, 
{\it Handbook of Brownian Motion}, 2nd~ed., 
Birkh{\" a}user, 2002.

\bibitem{BR} T.~Byczkowski and T.~Ryznar, 
Hitting distribution of geometric Brownian motion.
\emph{Studia Math.} \textbf{173}, (2006), 19--38.
%%MR{2204460}

\bibitem{HM-I} Y.~Hamana and H.~Matsumoto, 
The probability densities of the first hitting times 
of Bessel processes, 
J. Math-for-Industry, {\bf 4} (2012), 91--95.

\bibitem{HM-T} Y.~Hamana and H.~Matsumoto, 
The probability distributions of the first hitting times 
of Bessel processes, 
Trans. AMS, {\bf 365} (2013), 5237--5257.

\bibitem{HM-E} Y.~Hamana and H.~Matsumoto, 
Asymptotics of the probability distributions of 
the first hitting times of Bessel processes, 
Electron. Commun. Probab., {\bf 19} (2014), 1--5.  

\bibitem{K} J.~T.~Kent, 
Eigenvalue expansions for diffusion hitting times, 
Z. Wahr. Ver. Gebiete, {\bf 52} (1980), 309--319. 

\bibitem{MOS} W.~Magnus, F.~Oberhettinger and R.~P.~Soni, 
{\it Formulas and Theorems for the Special Functions of 
Mathematical Physics}, 3rd. ed., Springer--Verlag, 1966.

\bibitem{PY} J.~Pitman and M.~Yor, 
Bessel processes and infinitely divisible laws, 
in {\it Stochastic Integrals}, ed. by D.Williams, 285--370, 
Lecture Notes Math., vol. 851, Springer--Verlag, 1981.  

\bibitem{S} D.~Stroock, 
On the growth of stochastic integrals, 
Z. Wahr. Ver. Gebiete, {\bf 18} (1971), 340--344. 

\bibitem{U} K.~Uchiyama, 
Asymptotics of the densities of the first passage time distributions 
for Bessel diffusions, 
Trans. AMS, {\bf 367} (2015), 2719--2742.

\bibitem{W} G.~N.~Watson, 
{\it A Treatise on the Theory of Bessel Functions}, 
Reprinted of 2nd ed., Cambridge University Press, 1995.

\bibitem{wendel} J.~G.~Wendel, 
Hitting spheres with Brownian motion, 
Ann. Probab., {\bf 8} (1980), 164--169.  


\end{thebibliography}
\end{document}